\documentclass[a4paper,12 pt,leqno]{article}
\usepackage{amsmath,amsfonts,here,epsf}
\usepackage[latin1]{inputenc}
\usepackage[T1]{fontenc}
\begin{document}


\newtheorem{theo}{Theorem}[section]
\newtheorem{lemme}[theo]{Lemma}
\newtheorem{cor}[theo]{Corollary}
\newtheorem{problem}[theo]{Problem}
\newtheorem{prop}[theo]{Proposition}
\newtheorem{assu}[theo]{Assumption}
\newtheorem{nontheo}[theo]{Conjectured theorem}
\newtheorem{defi}[theo]{Definition}
\newtheorem{remark}[theo]{Remark}


\newcommand{\beq}{\begin{eqnarray}}
\newcommand{\enq}{\end{eqnarray}}
\newcommand{\be}{\begin{eqnarray*}}
\newcommand{\en}{\end{eqnarray*}}
\newcommand{\R}{\mathbb R}
\newcommand{\N}{\mathbb N}
\newcommand{\Linf}{L^{\infty}}
\newcommand{\dt}{\partial_t}
\newcommand{\Dt}{\frac{d}{dt}}
\newcommand{\demi}{\frac{1}{2}}
\newcommand{\ep}{^{\epsilon}}
\newcommand{\epu}{_{\epsilon}}
\newcommand{\vf}{\varphi}
\newcommand{\NN}{\mathbb N}
\newcommand{\RR}{\mathbb R}
\newcommand{\dx}{\partial_x}
\newcommand{\Sz}{{\mathcal{S}}}
\newcommand{\ds}{\displaystyle}
\newcommand{\fe}{f_\epsilon}
\newcommand{\T}{\cal T}
\newcommand{\U}{\cal U}
\let\cal=\mathcal
\newcommand{\rz}{\mathbb R}

\title{ \bf Contractive metrics for scalar conservation laws}
\date{}

\maketitle

\begin{center}Fran{\c c}ois Bolley\footnotemark[1], Yann Brenier\footnotemark[2],
  Gr\'egoire Loeper\footnotemark[3]\footnotemark[4]
\end{center}
\footnotetext[1]{Ecole normale sup\'erieure de Lyon, Umpa, 46, all\'ee d'Italie, F-69364 Lyon Cedex 07.
{\tt fbolley@umpa.ens-lyon.fr} }
\footnotetext[2]{Laboratoire J.A.Dieudonne,   UMR CNRS 6621,
           Universite de Nice Sophia-Antipolis,
           Parc Valrose,
           06108 Nice Cedex 2. {\tt brenier@math.unice.fr}}
\footnotetext[3]{Fields Institute \& University of Toronto}
\footnotetext[4]{Ecole Polytechnique F\'ed\'erale de  Lausanne, SB-IMA, 10015 Lausanne. {\tt gregoire.loeper@epfl.ch}}

\begin{abstract}
We consider nondecreasing entropy solutions to 1-d scalar conservation laws and show
that the spatial derivatives of such solutions satisfy a contraction
property with respect to the Wasserstein distance of any order. This result 
extends the  $L^1$-contraction property shown by Kru$\mathrm {\check z}$kov.
\end{abstract}

\vspace{8mm}

Existence and uniqueness of solutions to scalar conservation laws in one space
dimension have been established by 
Kru$\mathrm {\check z}$kov
in the framework of entropy solutions (see \cite{K} for instance),
 and among the properties satisfied by these solutions it is known
that the $L^1$ norm between any two of them is a non-increasing function of time.

In this work we shall focus on a class of entropy solutions such that a certain distance between
the space derivatives of any two such solutions is also nonincreasing in time.  On this class of solutions
this result extends the $L^1$ norm contraction property.
 
More precisely we consider as initial data nondecreasing functions on $\R$ with 
limits $0$ and $1$ at $-\infty$ and $+\infty$ respectively. These properties are
 preserved by the conservation law, and corresponding solutions have been shown in \cite{BG}
 to arise in some models of pressureless gases, obtained as a continuous limit of systems of sticky particles.
  Noticing that the distributional space derivative of these functions are
 probability
measures, we may consider the Wasserstein distance between the space derivatives of any two such
 solutions, and we shall prove in this paper that this distance is a nonincreasing 
 function of time, constant in the case of classical solutions.

\bigskip

\section{Introduction to the results}

\indent

Given a locally Lipschitz real-valued function $f$ on $\R$, called a flux,
we consider the scalar conservation law
\beq \label{cons law}
\left\{ \begin{array}{ll}
u_t + f(u)_x  \, = \, 0 \, , \qquad t > 0, \, \, x \in \R,   \\
u(0,.) \, =\,  u^{0},  
\end{array}
\right.
\enq
with unknown $u=u(t,x) \in \R$ and initial datum $u^{0} \in L^{\infty}(\R)$, and where the subscripts stand for derivation.

\medskip

We shall consider  solutions that are called {\it entropy solutions} (see \cite{S} for instance) and are
defined as follows:
a function $u=u(t,x) \in L^{\infty}( [0,+\infty[  \times \R )$ is said to be an entropy solution of (\ref 
{cons law})
on $  [0,+\infty[ \times \R$ if the entropy inequality 
\beq\label{entropy}
E(u)_t \, + \, F(u)_x \, \leq \, 0
\enq
holds in the sense of distributions for all convex Lipschitz function $E$ on $\R$, and with associated
flux  $F$  defined by
\begin{equation}\label{fluxF}
F(u) \, = \, \int_{0}^{u} f'(v) \, E'(v) \, dv.
\end{equation}
This means that
$$
\int_{0}^{+ \infty}  \! \! \int_{\R} \bigl( E(u) \, \varphi_t + F(u) \, \varphi_x  \bigr) \, dt\, dx \, + \, 
\int_{\R} E(u^0(x))\,  \varphi(0,x) \, dx \, \geq 0
$$
for all nonnegative $\varphi$ in the space ${\cal C}_{c}^{\infty}([0 , +\infty[ \times \R)$ of 
${\cal C}^{\infty}$ functions on $[0 , +\infty[ \times \R$ with compact support.

\medskip

We shall also consider {\it classical solutions}, that is, functions $u=u(t,x)$ in 
${\cal C}^{1}(]0,+\infty[ \times \R)  \cap \cal{C}([0,+\infty[ \times \R)$ satisfying
\eqref{cons law} pointwise.

In particular any classical solution to \eqref{cons law} satisfies \eqref{entropy}, i.e. is an
entropy solution,
and conversely any entropy solution satisfies \eqref {cons law} in the distribution sense.

\medskip

For entropy solutions, the following result is due to Kru$\mathrm {\check z}$kov (see \cite{K}):

\begin{theo} \label{Kruzkov}
For every $u^{0} \in L^{\infty}(\R)$, there exists a unique entropy solution $u$
to $(\ref{cons law})$ in 
$L^{\infty}([0,+\infty[ \times \rz) \cap \cal C ([0,+\infty[, L_{loc}^{1}(\R))$.
\end{theo}

Moreover for classical solutions, we have (see \cite{S} for instance):

\begin{theo} \label{Kruzkovclassique}
Given a ${\cal C}^{2}$ flux $f$ and a ${\cal C}^1$ bounded initial datum $u^0$ such that
$f' \circ u^0$ is nondecreasing on $\R$, the unique entropy solution $u$ to 
$(\ref{cons law})$ is a classical solution.
\end{theo}

In this work we shall consider initial data in the subset
${\cal U}$ of $L^{\infty}(\R)$  defined by

\begin{defi}
A function $v : \R \rightarrow \R$ belongs to ${\cal U}$ if it is nondecreasing, 
right-continuous,
and has limits $0$ and $1$ at $-\infty$ and  $+\infty$ respectively. 
\end{defi}

The following proposition expresses that this set is 
preserved by  the conservation law \eqref{cons law}:

\begin{prop}\label{Uconserve}
Given an initial datum  $u^0 \in {\cal U}$, the entropy solution $u$ given by 
Theorem  \ref{Kruzkov} is such that $u(t,.)$ belongs to
   ${\cal U}$ for all $t \geq 0$.
\end{prop}
   
More precisely, given any $t \geq  0$, the $L^{\infty}(\R)$ function $u(t,.)$ is $a.e.$ equal  to an element 
of the set ${\cal U}$, which on the other hand is characterized by

\begin{prop}\label{charact}
The distributional  derivative $v_x$ of any $v \in {\cal U}$ is a
Borel probability measure on $\R$, and  for any $x \in \R$,
$$
v(x) \, = \, v_x(]-\infty,x]).
$$
Conversely, if $\mu$ is a probability measure on $\R$, then $v$ defined on $\R$ as
$$
v(x) \, = \, \mu(]-\infty,x])
$$
belongs to ${\cal U}$, and $v_x = \mu$.
\end{prop}

Consequently the map $ v \mapsto v_x$ is  one-to-one from ${\cal U}$ onto the set ${\cal P}$ of 
probability measures on $\R$ (and ${\cal U}$ can be seen as the set of repartition functions 
of real-valued random variables).

Propositions \ref{Uconserve} and \ref{charact} allow us to characterize at any time the distance between
two solutions (with initial datum in ${\cal U}$) in terms of  their 
space derivatives, in particular by means of the Wasserstein distances: given
any real number $p \geq 1$, the Wasserstein distance of order $p$ is defined on the set of 
probability measures on $\R$ by 
$$
W_p(\mu,\tilde \mu)=\, \inf_{\pi} \, 
\Bigl( \,\int_{\R^2} \vert x-y \vert ^p \, d\pi(x,y)\,  \Bigr)^{1/p}
$$
where $\pi$ runs over the set of probability measures on $\R^2$ 
with marginals $\mu$ and
 $\tilde \mu$; these distances are considered here in a broad sense with possibly infinite values.

\medskip
 
This paper aims at proving that the Wasserstein distances between the space derivatives of any two such entropy
solutions is a nonincreasing function of time:

\begin{theo} \label{theoreme}
Given  a locally Lipschitz real-valued function $f$ on $\R$ and two initial data $u^0$ and ${\tilde u}^0$ 
in ${\cal U}$, let $u$ and 
${\tilde u}$ be the associated entropy solutions to \eqref{cons law}. Then, for any 
 $t \geq 0$ and $p \geq 1$, we have (with possibly infinite values)
$$
 W_p(u_{x}(t,.), {\tilde u}_{x}(t,.)) \, \leq \, W_p(u_{x}^{0}, {\tilde u}_{x}^{0}).
 $$
\end{theo}

We shall see in Section \ref{section wasserstein} that for $p=1$ the distance $W_1$ satisfies
$$
W_1(v_x,{\tilde v}_x) \, = \, \Vert v-{\tilde v} \Vert_{L^1(\R)}
$$
for all $v,{\tilde v} \in {\cal U}$. Hence Theorem \ref{theoreme} reads in the case $p=1$:
$$
\Vert u(t,.) - {\tilde u}(t,.) \Vert_{L^1(\R)} \, \leq \, \Vert u^0 - {\tilde u}^0 \Vert_{L^1(\R)}.
$$
Thus, for initial profiles in ${\cal U}$, we recover the $L^1$-contraction property given by
Kru$ {\check z}$kov.

\medskip

In the case of classical solutions, the result of Theorem \ref{theoreme} is improved, since the Wasserstein
 distance between two solutions is conserved: 
\begin{theo}\label{theoremeclassique}
Given  a ${\cal C}^1$ real-valued function $f$ on $\R$, let $u^0$ and ${\tilde u}^0$ in ${\cal U}$ be
 two initial data such that  the associated entropy solutions 
$u$ and ${\tilde u}$  to \eqref{cons law} are classical solutions,  increasing in $x$ for all $t\geq 0$.
 Then for any  $t \geq 0$ and $p \geq 1$ we have (with possibly infinite values)
$$
W_p(u_{x}(t,.), {\tilde u}_{x}(t,.)) \, = \, W_p(u_{x}^{0}, {\tilde u}_{x}^{0}).
$$
\end{theo}

\medskip

>From these general results can be induced some corollaries in the case of initial data in the subsets
${\cal U}_p$ of ${\cal U}$ defined as:
\begin{defi}
Let $p \geq 1$. A function $v$ in ${\cal U}$ belongs to ${\cal U}_p$ if its distributional
derivative $v_x$ has finite moment of order $p$, that is, if $\displaystyle \int_{\R} \vert x \vert^p \, dv_x(x)$ 
is finite.
\end{defi}

As in Proposition \ref{charact} the map $ v \mapsto v_x$ is  one-to-one from ${\cal U}_p$ onto the set 
${\cal P}_p$ of probability measures on $\R$ with finite moment of order $p$. But we shall note in 
Section \ref{section wasserstein} that the map $W_p$ on ${\cal P}_p \times {\cal P}_p$
defines a distance on ${\cal P}_p$.
Then the real-valued map $d_p$ defined on ${\cal U}_p \times {\cal U}_p$ by
$$
d_p(v, {\tilde v}) \, = \, W_p(v_x, {\tilde v}_x)
$$
induces a distance on ${\cal U}_p$, and for the associated topology we have
\begin{cor}\label{corsolutionUp}
Given a locally Lipschitz function $f$ on $\R$, $p \geq 1$ and $u^0 \in {\cal U}_p$,  the entropy solution $u$ to 
\eqref{cons law} belongs to ${\cal C}([0,+\infty[,{\cal U}_p).$
\end{cor}
In particular for $p=1$
$$
d_1(v,{\tilde v}) = W_1(v_x,{\tilde v}_x) = \Vert v- {\tilde v} \Vert_{L^1(\R)},
$$
and the previous result can be precised by
\begin{cor}\label{corsolutionU1}
Given a locally Lipschitz function $f$ on $\R$ and $u^0 \in {\cal U}_1$, the entropy solution $u$ to  \eqref{cons law} is such that
$$
\Vert u(t,.) - u(s,.) \Vert_{L^1(\R)} \leq \vert t-s \vert \, \Vert f' \Vert_{L^{\infty}([0,1])}.
$$
\end{cor}
This known result holds under weaker assumptions (for $u^0$ with bounded variation, see \cite{S}), but in
our case it will be recovered in a straightforward way.

\medskip

Finally Theorem \ref{theoremeclassique} can be precised in the ${\cal U}_p$ framework in 
the following way:
\begin{cor}\label{corclassique}
Given a ${\cal C}^2$ convex flux $f$ and two ${\cal C}^1$ increasing initial data $u^0$ and ${\tilde u}^0$ 
in ${\cal U}_p$
for some $p \geq 1$, the following three properties hold:

1. the associated entropy solutions $u$ and ${\tilde u}$ are classical solutions;

2.  $u(t,.)$ and ${\tilde u}(t,.)$ belong to ${\cal U}_p$ and are increasing for all $t \geq 0$;

3. for all $t \geq 0$, we have (with finite values)
$$
W_p(u_{x}(t,.), {\tilde u}_{x}(t,.)) \, = \, W_p(u_{x}^{0}, {\tilde u}_{x}^{0}).
$$
\end{cor}

\medskip
The paper is organized as follows.
 The definition and some properties of Wasserstein distances are discussed in greater detail in Section 
\ref{section wasserstein}. In Section \ref{sectionpreuveclassique} we consider the case of
classical solutions, proving Theorem \ref{theoremeclassique} and
Corollary \ref{corclassique}. Then the general case of entropy solutions is  studied 
in Sections \ref{sectiondiscretisation} and \ref{sectionpreuve}:  
more precisely  in Section \ref{sectiondiscretisation} we  introduce a time-discretized scheme, show
the $W_p$ contraction  property for this discretized evolution and prove the convergence of the corresponding approximate solution toward the 
entropy solution;
Theorem \ref{theoreme}  and its corollaries follow from this in Section \ref{sectionpreuve}. In Section
 \ref{sectionviscous} we shall finally see how such results extend to viscous conservation laws.

\section {Wasserstein distances}\label{section wasserstein}
\indent

In this section $p$ is a real number with $p\geq 1$, ${\cal P}$ (resp. ${\cal P}_p$) stands for the
set of probability measures on $\R$ (resp. with finite moment of order $p$) and $dx$ 
for the Lebesgue measure on $\R$.
 
The Wasserstein distance of order $p$, valued in $\R \cup \{+\infty\}$, is  defined on 
${\cal P} \times {\cal P}$ by
\beq\label{defwas1}
W_p(\mu,{\tilde \mu}) \, = \, \inf_{\pi} \, 
\Bigl( \,\int_{\R^2} \vert x-y \vert ^p \, d\pi(x,y)\,  \Bigr)^{1/p}
\enq
where $\pi$ runs over the set of probability measures on $\R^2$ with marginals 
$\mu$ and ${\tilde \mu}$. It is equivalently defined   by
\beq\label{defwas2}
W_p(\mu, {\tilde \mu})= 
\inf_{X_{\mu},X_{{\tilde \mu}}}
\Bigl(\int_{0}^{1}\left|X_{\mu}(w)-X_{{\tilde \mu}}(w)\right|^p \, dw
\Bigr)^{1/p}
\enq
where the infimum is taken over all random variables $X_{\mu}$ and $X_{{\tilde \mu}}$ on 
the probability space $(]0,1[,dw)$
with respective laws  $\mu$ and ${\tilde \mu}$.
It takes finite values on ${\cal P}_p \times {\cal P}_p$
and indeed defines a distance on  ${\cal P}_p$.

For complete references about the Wasserstein distances and related topics the reader
 can refer to  \cite{V}. We only mention that both infima in (\ref{defwas1}) and (\ref{defwas2})
  are achieved, and 
 for the second definition we shall precise some random variables that achieve the infimum. For this
 purpose we introduce the notion of generalized inverse:

\begin{defi}\label{ginverse}
Let $v$  belong to ${\cal  U}$. Then its generalized inverse is the function $v^{-1}$
defined on $]0,1[$ by 
\be
v^{-1}(w) \, = \, \inf \{x \in \R; \, v(x) > w\}.
\en 
\end{defi}

Then $v^{-1}$ is a nondecreasing random variable on $(]0,1[, dw)$ by definition,
with law $v_x$ since
\be
\int_{0}^{1} f(v^{-1}(w)) \ dw 
&=&  
 \int_{0}^{1}  \Bigl( \int_{\R} f'(s) {\bf 1}_{\{s \leq v^{-1}(w)\}} \, ds  \Bigr) \, dw\\
&=&
\int_{\R} \Bigl( \int_{0}^{1} f'(s) {\bf 1}_{\{v(s) \leq w\}} \, dw  \Bigr) \, ds \\
&=& 
\int_{\R} f'(s)(1-v(s)) \ ds \\
&=& \int_{\R} f(s) \ dv_x(s)  
\en
for all $f$ in ${\cal C}^{1}_{c}(\R)$. In particular its repartition function is $v$.

Moreover this generalized inverse achieves the infimum in (\ref{defwas2}):
\begin{prop}\label{propertyW_p}
Let $v$ and ${\tilde v}$ in ${\cal U}$. Then we have
 (with possibly infinite values)
\be 
W_p(v_x,{\tilde v}_x) \, = \, \Bigl( \int_{0}^{1} \vert v^{-1}(w) - 
       {\tilde  v}^{-1}(w) \vert^p \, dw \Bigr)^{1/p}
\en
for all $p \geq 1$. In  particular for $p=1$ we also have
\be
W_1(v_x,{\tilde v}_x) \, = \, \Vert v - {\tilde v} \Vert_{L^1(\R)}.
\en 
\end{prop}

\medskip
{\it Proof.}
The general result is proved in \cite{V}.  The result specific to the case $p=1$ follows 
by introducing, for a given $v \in {\cal U}$, the map defined on $\R \times ]0,1[$ by
$$
jv (x,w) = \left\{ \begin{array}{ll}
                           1 & \mbox{if $v(x) > w$}\\
                           0 & \mbox{if $v(x) \leq w,$}   
                          \end{array}
                  \right.
$$
for which we have
$$
|v^{-1}-\tilde v^{-1}|(w)  = \int_{\R}|jv - j\tilde v|(x,w) \  dx
$$
for almost every $ w \in \, ]0,1[$, and
$$
\int_{0}^{1}|jv - j\tilde v|(x,w) \ dw\, = \, \vert v-\tilde v \vert(x)
$$
for almost every $x \in \R$. Integrating the first equality on $w$ in $]0,1[$ and the second one on $x$ 
in $\R$, we deduce
$$
\int_{0}^{1}|v^{-1}-\tilde v^{-1}|(w) \ dw = \int_{\R}|v-\tilde v|(x) \ dx.  
$$
$\hfill\Box$

\bigskip

Given $v \in {\cal U}$, its generalized inverse $v^{-1}$ is actually the a.e. unique nondecreasing 
random variable on $(]0,1[, dw)$ with law $v_x$. Given any other random variable $X$ on $(]0,1[, dw)$
with law $v_x$, $v^{-1}$ is called the (a.e. unique) nondecreasing rearrangement of $X$ (see \cite{V}).

\medskip

We conclude this section recalling a result relative to the convergence of probability measures.
 A sequence $(\mu_n)$ of probability measures on $\R$ is said to converge weakly toward 
a probability measure $\mu$ if, as $n$ goes to $+\infty$, 
$\displaystyle \int_{\R} \varphi \, d\mu_n$ tends to 
$\displaystyle \int_{\R} \varphi \, d\mu$  for all bounded continuous real-valued  functions
 $\varphi$ on $\R$ (or equivalently for all ${\cal C}^{\infty}$ functions $\varphi$ with compact
 support, that is, if $\mu_n$ converges to $\mu$ in the distribution sense). Given $p\geq 1$ this convergence
is metrized on ${\cal P}_p$ by the distance $W_p$ as shown by the following proposition (see \cite{V}):

\begin{prop}\label{cvWp/weak}
Let $p \geq 1$, $(\mu_n)$ a sequence of probability measures in ${\cal P}_p$ and $\mu \in {\cal P}.$
Then the following statements are equivalent:

i) $(W_p(\mu_n,\mu))$ converges to $0$;

ii) $(\mu_n)$ converges weakly to $\mu$ and $\displaystyle \sup_{n} \int_{\vert x \vert \geq R}
 \vert x \vert^p \, d\mu_n(x)$ tends to $0$ as $R$ goes to infinity.
\end{prop}

In this proposition we do not a priori assume that $\mu$ belongs to ${\cal P}_p$, but it can be noted that
this property is actually induced by  any of both hypotheses $i)$ and $ii)$. 

For measures in $\cal P$ we have the  weaker result:

\begin{prop}\label{weakcv}
Let $p \geq 1$, $(\mu_n)$ and $(\nu_n)$ two sequences in ${\cal P}$ 
converging weakly to $\mu$ and $\nu$ in ${\cal P}$ respectively. Then 
(with possibly infinite values)
$$
W_p(\mu,\nu) \leq \liminf_{n \to +\infty} W_p(\mu_n,\nu_n).
$$
\end{prop}

\section{The  case of classical solutions:  Theorem \ref{theoremeclassique}
and corollary}\label{sectionpreuveclassique}
\indent

\subsection{Proof of Theorem \ref{theoremeclassique}}
\indent 

We consider two 
classical solutions $u$ and ${\tilde u}$ to (\ref{cons law}) such that
 $u(t,.)$ and ${\tilde u}(t,.)$ belong to ${\cal U}$ and are increasing for all $t \geq 0$, and we shall
prove that
$$
W_p(u_{x}(t,.), {\tilde u}_{x}(t,.)) \, = \, W_p(u_{x}^{0}, {\tilde u}_{x}^{0})
$$
as a consequence of Proposition \ref{propertyW_p}.

The map $u^0$ is increasing from $0$ to $1$, so has a (true) inverse 
$X(0,.)$ defined on $]0,1[$ by
$$
u^0(X(0,w)) \, = \, w. 
$$
Then, given $w \in \, ]0,1[$, we consider a characteristic curve $t \mapsto X(t,w)$ solution of
\begin{equation}\label{caracteristiques}
X_t(t,w) = f'\bigl(u(t,X(t,w))\bigr)
\end{equation}
for $t \geq 0$, and taking value $X(0,w)$ at $t=0$. 
Since $f$ is ${\cal C}^1$ and $u$ is bounded there exists a (non necessarily unique) solution 
$X(.,w)$ to 
\eqref{caracteristiques} by Peano Theorem (see \cite{H} for instance); 
moreover by a classical
computation from \eqref{cons law} it is known to satisfy
\begin{equation}\label{inverse}
u(t,X(t,w)) \, = \, w
\end{equation}
for all $t \geq 0$, from which it follows that
$$
X_t(t,w) \, \Bigl( = f'(u(t,X(t,w)))  \Bigr)= f'(w)
$$
and hence
\begin{equation}\label{X(t)classique}
X(t,w) \, = \, X(0,w) \, + \, tf'(w).
\end{equation}
In particular there exists  a unique
 solution $X(.,w)$ to \eqref{caracteristiques}. Now given $t \geq 0$, $X(t,.)$ is the (true) inverse of the increasing function $u(t,.)$
(by \eqref{inverse}), and Proposition \ref{propertyW_p} writes
$$
W_{p}(u_{x}(t,.),{\tilde u}_{x}(t,.)) \, = \, \Bigl( \int_{0}^{1} \vert X(t,w) - 
       {\tilde  X}(t,w) \vert^{p} \, dw \Bigr)^{1/p}.
$$
But from \eqref{X(t)classique} we obtain
\begin{equation}\label{cst}
X(t,w) \, - \, {\tilde X}(t,w) \, = \, X(0,w) \, - \, {\tilde X}(0,w).
\end{equation}
This result ensures in particular that $\ds W_p(u_{x}(t,.),{\tilde u}_{x}(t,.))$  
remains constant in time, may its initial value be finite or not; note however
that \eqref{cst} is actually much stronger that 
Theorem \ref{theoremeclassique}.

\bigskip

\subsection{Proof of Corollary \ref{corclassique}}
\indent

We assume that $f$ is a ${\cal C}^2$ convex function on $\R$, and $u^0$ is
a ${\cal C}^1$ increasing initial profile in ${\cal U}_p$.

First of all we note that the associated entropy solutions $u$ is a classical solution
in view of Theorem \ref{Kruzkovclassique}: this result is proved in \cite{S} for instance,
and its proof also ensures that  $u(t,.)$ is increasing for all $t \geq 0$.

Then we check that the moment property is preserved by the conservation law, that is, that
$u(t,.)$ also  belongs to ${\cal U}_p$ for any $t \geq 0$. Indeed, given $t \geq 0$, we have by the
change of variable $w = [u(t,.)](x)$:
\begin{eqnarray*}
\int_{\R} \vert x \vert^p \, u_{x}(t,x) \, dx
& = & \int_{0}^{1} \vert X(t,w) \vert^p \, dw \\
& = &  \int_{0}^{1} \vert X(0,w) + t f'(w) \vert^p \, dw \\
& \leq & 2^{p-1} \left[ \int_{0}^{1} \vert X(0,w) \vert^p \, dw \, +
 t^p  \Vert f' \Vert^p_{L^{\infty}(]0,1[)}
\right]
\end{eqnarray*}
which is finite since
$$
\int_{0}^{1} \vert X(0,w)\vert^p \, dw = \int_{\R} \vert x \vert^p \, u_{x}^{0}(x) \, dx
$$
is finite by assumption.
This ends the proof of Corollary \ref{corclassique}. $\hfill\Box$

\section{Time discretization of the conservation law}\label{sectiondiscretisation}
\indent

 In the  previous section we have seen that the classical
   solutions are obtained through the  method of characteristics, that we now summarize in our case:
    given an
  initial profile $u^0$ in ${\cal U}$ such that the corresponding solution $u$ is ${\cal C}^{1}$ and 
  increasing 
 in $x$ for all $t \geq 0$, let $X(0,.)$ be its inverse, defined by
 $$
 u^0(X(0,w)) \, = \, w
 $$
 for all $w \in \, ]0,1[$. Let then  $X(0,w)$ evolve into
\begin{equation}\label{X(t)}
 X(t,w) \, = \, X(0,w) + t f'(w)
\end{equation}
for all $ t \geq 0$ and $w \in \, ]0,1[$ (see  \eqref{X(t)classique}). 
 The solution $u(t,.)$ is then the inverse of the increasing map $X(t,.)$, that is, is 
 the unique solution of
 $$
 X(t,u(t,x)) \, = \, x.
 $$

\medskip

In the general case, defining $X(0,.)$ in some similar way, there is no hope for
 the function $X(t,.)$
defined by \eqref{X(t)}
to be increasing for $t > 0$; inverting it would thus lead to a multivalued function, and 
no more to the
entropy solution of the conservation law, as in the particular case discussed above.

However, averaging (or "collapsing") this multivalued function into a single-valued function, Y. Brenier showed in 
\cite{B} how to build an approximate solution to the conservation law.

We now precisely describe this so-called Transport-Collapse method in our case.

 \subsection{Definition and $W_p$ contraction property of the discretized solution}\label{Time discretization 
 of the problem}
\indent

Let  $u^{0} \in {\cal U}$ be some fixed initial profile, with generalized inverse
 $X(0,.)$ given as in Definition \ref{ginverse} by
 $$
 X(0,w) \, = \, \inf \{x \in \R; u^{0}(x) > w \}
 $$
 for all $w \in \, ]0,1[$.  $X(0,.)$ can be seen as a random variable on the
 probability space $]0,1[$ equipped with the Lebesgue measure $dw$; its law is
 $u^{0}_x$, as pointed out after Definition \ref{ginverse}.
  
 We let then $X(0,.)$ evolve according to the method of
 characteristics, denoting
 $$
 X(h,w) \, = \, X(0,w) \, + \, hf'(w)
 $$
 for all $h \geq 0$ and almost every $w \in \, ]0,1[$. Again, given $h \geq 0$, $X(h,.)$ can be seen as a random variable on 
 $]0,1[$; let then
 $T_h u^0$ be its repartition function, that is, the function belonging to ${\cal U}$ and defined at any $x \in \R$
 as the Lebesgue measure of the set $\{w\in \, ]0,1[; \, X(h,w) \leq x \}$.
It is given by
 $$
 T_h u^0(x) \, = \,  \int_{0}^{1} {\bf 1}_{\{X(h,w) \leq x\}}(w) \, dw.
 $$

\medskip

We summarize  this construction in the following definition:
\begin{defi}\label{defiTh}
Let $v\in{\cal U}$ with generalized inverse $X(0,.)$ defined on $]0,1[$ by
$$
X(0,w) \, = \, \inf \{x \in \R; v(x) > w \}.
$$
Then, given $ h \geq 0$, and letting 
$$
X(h,w) \, = \, X(0,w) \, + \, hf'(w)
$$ 
for almost every $w \in \, ]0,1[$, we define the ${\cal U}$ function $T_hv$ on $\R$ by
$$
 T_h v(x) \, = \, \int_{0}^{1} {\bf 1}_{\{X(h,w) \leq x\}}(w) \, dw.
$$ 
\end{defi}

\medskip

In the case of Section \ref{sectionpreuveclassique} (see  \eqref{X(t)classique}),
 it turns out that $X(h,.)$ is the (true) inverse of $T_hu^0,$ and $(h,x) \mapsto T_hu^{0}(x)$ is 
 exactly
  the entropy solution to equation $(\ref {cons law})$ with initial datum $u^{0}$ in $\cal U$.
 This does not hold anymore in the general case, 
but will allow us to build an approximate solution $S_hu^0$  by iterating the operator  $T_h$.  Let us first give two important properties of $T_h$:

\begin{prop}\label{propTh}
Let $h\geq 0$, $T_h$ defined as above and $p \geq 1$. Then 

i) $T_hv$ belongs to ${\cal U}_p$  if so does  $v$.

ii) For any $v$ and ${\tilde v}$ in ${\cal U}$ we have (with possibly infinite values unless
$v$ and ${\tilde v} \in {\cal U}_p$)
$$ 
 W_p([T_h v ]_x, [T_h {\tilde v} ]_x) \, \leq 
 \, W_p(v_x, {\tilde v}_x).
 $$
 \end{prop}

\medskip

{\it Proof.} It is  really similar to what has been done in  Section
\ref{sectionpreuveclassique} 
as for Corollary \ref{corclassique}.

$i)$ $T_h v$  belongs to ${\cal U}$ as a repartition function of a random variable, and   we have
\begin{eqnarray*}
\int_{\R} \vert x \vert^p \, d[T_hv]_x(x) 
& =  & \int_{0}^{1} \vert X(h,w) \vert^p \, dw \\
&= &\int_{0}^{1} \vert X(0,w) \, + \, hf'(w) \vert^p \, dw \\ 
&\leq & 2^{p-1} \int_{0}^{1} \vert X(0,w) \vert^p +   \vert h f'(w) \vert^p \, dw \\
&\leq & 2^{p-1} \left[ \int_{\R} \vert x \vert^p \, dv_{x}(x) +
 h^p \Vert f' \Vert^p_{L^{\infty}(]0,1[)}
\right],
\end{eqnarray*}
which ensures that $[T_hv]_x$ has finite moment of order $p$ if so does $v_{x}$.

\bigskip

$ii)$  On  one hand the generalized inverses $X(0,.)$
 and $\tilde{X}(0,.)$ of  $v$
 and ${\tilde v}$ respectively satisfy
\begin{equation}\label{time0}
 W_{p}(v_x, {\tilde v}_x) \, = \, \Bigl( \int_{0}^{1} \vert X(0,w) - \tilde{X}(0,w) \vert^p 
 \, dw \Bigr)^{1/p}
\end{equation}
 by Proposition \ref {propertyW_p} (with finite values if both $v$ and ${\tilde v}$ belong to
 ${\cal U}_p$, and possibly infinite  otherwise).
On the other hand 
 $X(h,.)$ and $\tilde{X}(h,.)$ have respective 
 law $[T_h v]_x$ and $[T_h {\tilde v}]_x$, so
\begin{equation}\label{timet}
W_{p}([T_h u^{0}]_x, [T_h {\tilde u}^{0}]_x) \, 
\leq  \, \Bigl(\int_{0}^{1} \vert X(h,w) - \tilde{X}(h,w) \vert^p  \, dw \Bigr)^{1/p}
\end{equation}
by definition of the Wasserstein distance. 
But 
$$
 X(h,w) - \, \tilde{X}(h,w) \, = \, X(0,w) \, - \, \tilde{X}(0,w)
 $$
 for almost every $w \in \, ]0,1[$ by definition, which concludes the argument by 
 \eqref{time0} and \eqref{timet}.
 
Note again that \eqref{timet}  holds only as an inequality since $X(h,.)$ and ${\tilde X}(h,.)$ are
not necessarily nondecreasing, which was the case in the example
discussed in Section \ref{sectionpreuveclassique}.
$\hfill\Box$

\medskip

We now use the operator $T_h$ defined above to build an approximate solution $S_hu^0$
to the conservation law
$(\ref {cons law})$:

\begin{defi}\label{defiSh}
Let  $h$ be some  positive number and $v \in {\cal U}$.
For any $t \geq 0$ decomposed as $t = (N+s)h$ with  $N \in \N$ and $0 \leq s < 1$,
we let
$$
S_hv(t,.) \, = \, (1-s) \, T_{h}^{N} v(.) \, + \, s \, T_{h}^{N+1} v(.)
$$
where $T_{h}^{0}v = v$ and $T_{h}^{N+1}v = T_h(T_{h}^{N}v)$.
\end{defi}

These iterations make sense because $T_hv\in{\cal U}$ if $v\in {\cal U}$,
$S_hv(t,.) \in {\cal U}$ (resp. ${\cal U}_p$) for any $h, t \geq 0$ and
 $v \in {\cal U}$ (resp. ${\cal U}_p$). 

We now prove two contractions properties on these approximate solutions. We first have the $L^1(\R)$ 
contraction property:
 
\begin{prop}\label{contractL1}
Let $h$ be some fixed positive number and $S_h$ defined as above. Then, for any $v \in {\cal U}$ and
$s,t \geq 0$ we have
$$
\Vert S_hv(t,.) - S_hv(s,.) \Vert_{L^1(\R)} \, \leq \, \vert t-s \vert \, \Vert f' \Vert_{L^{\infty}([0,1])}.
$$
\end{prop}

\medskip
{\it Proof.}
As in  \cite{B} we first observe that $\|T_hV -V\|_{L^1(\R)}\leq h \Vert f' \Vert_{L^{\infty}(]0,1[)}$
 for any $V \in {\cal U}$, then let 
 $t=(M+\mu)h$ and $s= (N+\nu)h$  with $M, N \in \N$ and $0 \leq \mu, \nu <1$, 
and, for instance assuming that $M>N$, prove the proposition by applying this first bound
to $V = S_hv(kh,.) = T_{h}^{k}v$ for $k=N+1, \dots, M-1$.  
$\hfill \Box$

\medskip

Then we have the $W_p$ contraction property:

\begin{prop}\label{Shcontract}
Let $h$ be some fixed positive number and $S_h$ defined as above. Then, 
given $v$ and ${\tilde v}$ in 
${\cal U}$, we have for any $t \geq 0$:
$$
W_p([S_h v]_x(t,.), [S_h{\tilde v}]_x(t,.)) \, \leq 
 \, W_p(v_x, {\tilde v}_x).
 $$
\end{prop}

\medskip

{\it Proof.}
It follows from  Proposition \ref{propTh} (about $T_h$) and to the convexity of the $W_p$ distance
 to the power $p$, in the sense that
 $$
 W_{p}^{p}(\alpha \mu_1 + (1-\alpha) \mu_2, \alpha \nu_1 + (1-\alpha) \nu_2) \, \leq \, 
 \alpha W_{p}^{p}(\mu_1 , \nu_1 ) + ( 1-\alpha) W_{p}^{p}(\mu_2, \nu_2)
 $$
 for all real number $\alpha \in [0,1]$ and probability measures $\mu_1, \mu_2, \nu_1$ and $\nu_2$
  (see \cite{V} for instance). 
$\hfill\Box$

\bigskip

We shall now recall the convergence of the scheme toward the entropy solution of the conservation
law.

\subsection{Convergence of the scheme in the $L^{1}_{loc}(\R)$ sense}\label{cvL1loc}
\indent

In this section we consider the space $\cal C([0,+\infty[, L^1_{loc}(\R))$ equipped with the topology
defined by the semi-norms 
$$
q_{nm}(f) \, = \, \sup_{t \in [0,n]} \, \int_{-m}^{m} \vert f(x) \vert \, dx
$$
for any integers $n$ and $m$ and $f \in {\cal C}([0,+\infty[, L^1_{loc}(\R))$. Then we have
\begin{prop} \label{schemeconvergence}
 Let $u^0 \in {\cal U}$. Then, as $h$ goes to $0$, the function
$S_{h}u^0$ converges in ${\cal C}([0,+\infty[, L^1_{loc}(\R))$ to the entropy solution 
of (\ref{cons law}) with initial datum $u^0$.
\end{prop} 

We briefly give the steps of the proof, which follows the one of Brenier in \cite{B}, adapted  
to functions of ${\cal U}$ instead of $L^1(\R)$.

We first prove that the family $(S_hu^0)_h$ is 
relatively compact in ${\cal C}([0,+\infty[, L^1_{loc}(\R))$ by means of
Proposition \ref{contractL1}, and Helly and Ascoli-Arzela Theorems. Then we check that the
 limit of any sequence of $(S_hu^0)_h$ converging in 
${\cal C}([0,+\infty[, L^1_{loc}(\R))$ is an entropy solution to the conservation law
\eqref{cons law} with initial datum $u^0$. By the uniqueness of this solution
ensured by Theorem \ref{Kruzkov}, this concludes the proof of Proposition \ref{schemeconvergence}.

\subsection{Convergence of the scheme in $W_p$ distance sense}\label{subsectionpreuveconvergenceWp}
\indent

We first prove a uniform equiintegrability result on the approximate solutions:

\begin{prop}\label{equiintegrable}
Let $S_h$ be defined as above,  $v\in  {\cal U}_p$ and $T \geq 0$. Then 
\be
\sup_{0 \leq h \leq T} \sup_{0 \leq t \leq T} \int_{|x|\geq R} |x|^p  d[S_hv]_x(t,x) 
\en
tends to $0$ as $R$ goes to infinity.
\end{prop}

{\it Proof.}
We again denote $\ds M = \Vert f' \Vert_{L^{\infty}(]0,1[)}$, and first consider $T_h$ itself, writing
\begin{eqnarray*}
\int_{\vert x\vert \geq R} \vert x \vert^p d[T_hv]_x(x) 
& =& \int_{0}^{1} \vert v^{-1}(w) +hf'(w) \vert^p \, {\bf 1}_{\{\vert v^{-1}(w) +hf'(w) \vert \geq R\}} \, dw \\
& \leq & \int_{\R} (\vert x \vert +hM)^p {\bf 1}_{\{\vert x \vert +hM \geq R\}} \, dv_x(x)  \\
& \leq & \left( 1+\frac{hM}{R-hM} \right)^p \int_{\vert x \vert \geq R-hM} \vert x \vert^p \, dv_x(x)
\end{eqnarray*}
for $R  >hM$. From this computation we deduce by iteration
$$
\int_{\vert x\vert \geq R} \vert x \vert^p d[T_{h}^{N}v]_x(x) 
 \leq \, \prod_{j=1}^{N} \left( 1 + \frac{hM}{R-jhM} \right)^p \int_{\vert x \vert \geq R-NhM}
 \vert x \vert^p \, dv_x(x)
$$
for $R > NhM$, with
$$
\prod_{j=1}^{N} \left( 1 + \frac{hM}{R-jhM} \right) \leq \left( 1 + \frac{hM}{R-NhM} \right)^{N} \leq
\exp \left( \frac{NhM}{R-NhM} \right).
$$

Thus
$$
\int_{\vert x\vert \geq R} \vert x \vert^p d[S_{h}v]_x(Nh,x) 
 \leq \exp \left( \frac{pTM}{R-TM} \right) \int_{\vert x \vert \geq R-TM}
 \vert x \vert^p \, dv_x(x)
$$
for any $N$ and $h$ such that $Nh \leq T$.

>From this we get for instance
$$
\int_{\vert x\vert \geq R} \vert x \vert^p d[S_{h}v]_x(t,x)  
 \leq \exp \left( \frac{2pTM}{R-2TM} \right) \int_{\vert x \vert \geq R-2TM}
 \vert x \vert^p \, dv_x(x)
$$
for any $t$ and $h$ smaller than $T$. This concludes the argument since the last integral tends to $0$ as 
$R$ goes to infinity.
$\hfill\Box$

\bigskip

>From this we deduce the convergence of the scheme in $W_p$ distance sense:

\begin{prop}\label{convergenceWp}
Let $u^0 \in {\cal U}_p$  and
 $u$ be the entropy solution 
to (\ref{cons law}) with initial datum $u^0$.
 Then, for any $t \geq 0$, 
$W_p([S_h u^0]_x(t,.), u_x(t,.) )$ converges to $0$ as  $h$ goes to $0$.
\end{prop}

\medskip

{\it Proof.}
Given $t \geq 0$, $S_hu^0(t,.)$ converges to $u(t,.)$ in $L^{1}_{loc}(\R)$ as $h$ goes to $0$ (by Proposition 
\ref{schemeconvergence}), so $[S_hu]_x(t,.)$ converges to the probability measure $u_x(t,.)$, first
in the distribution sense, then in the weak sense of probability measures, and finally in $W_p$ distance by
 Propositions \ref{equiintegrable}  and \ref{cvWp/weak}.

Note in particular that  $u_x(t,.)$ has finite moment of order $p$ for any $ t \geq 0$, that is, 
$u(t,.)$ belongs to ${\cal U}_p$.
$\hfill\Box$

\section{The general case of entropy solutions: Theorem \ref {theoreme} and corollaries}\label{sectionpreuve}

\subsection{Proof of Theorem \ref {theoreme}}\label{soussectionpreuve}

\indent

 We let $p \geq 1$ and consider two initial data $u^0$ and ${\tilde u}^0$ in ${\cal U}$ with
 associated entropy solutions $u$ and $\tilde u$.

 Given $t \geq 0$,   Proposition
\ref{schemeconvergence} yields again the convergence of $[S_h u^0]_x(t,.)$  to $u_x(t,.)$ 
 in the weak sense of probability measures. Since this holds also
for $\tilde u^0$, we obtain
$$
W_p(u_x(t,.), {\tilde u}_x(t,.)) \leq \liminf_{h \to 0}
 W_p([S_h u^0]_x(t,.), [S_h {\tilde u}^0]_x(t,.))
$$
by Proposition \ref{weakcv}. But, for each $h$,
$$
W_p([S_h u^0]_x(t,.), [S_h {\tilde u}^0]_x(t,.)) \leq W_p(u_x^0, {\tilde u}_x^0)
$$
by Proposition \ref{Shcontract}, so finally
$$
W_p(u_x(t,.), {\tilde u}_x(t,.)) \leq W_p(u_x^0, {\tilde u}_x^0).
$$
This concludes the argument.

\subsection{Proof of Corollary \ref{corsolutionUp}}

\indent

We recall that in the introduction we have defined a distance on each ${\cal U}_p$ by letting
$$
d_p(u, {\tilde u}) \, = \, W_p(u_x, {\tilde u}_x),
$$
and we now prove that, given $p \geq 1$ and $u^0 \in {\cal U}_p$, the entropy solution $u$ to the
conservation law \eqref{cons law} belongs to ${\cal C}([0,+\infty[, {\cal U}_p)$.

We first note, in view of the proof of Proposition \ref{convergenceWp}, that $u(t,.)$ indeed belongs to
${\cal U}_p$ for all $t \geq 0$.

Then, given $s \geq 0$, we need to prove that $d_p(u(t,.),u(s,.))\, (=W_p(u_x(t,.),u_x(s,.)))$ tends to $0$
as $t$ goes to $s$. Indeed, on one hand $u(t,.)$ tends to $u(s,.)$ in $L^{1}_{loc}(\R)$ by Theorem \ref{Kruzkov}, so $u_x(t,.)$ tends to
$u_x(s,.)$, first in the distribution sense, then in the weak sense of probability measures.

On the other hand, given $T > s$, we now prove that 
$\ds \sup_{0 \leq t \leq T} \int_{\vert x \vert \geq R} \vert x \vert^p \, du_x(t,x)$ goes to $0$ as $R$ goes
 to infinity. For this, given $\varepsilon > 0$, let $R$ such that
$$
\sup_{0 \leq h \leq T} \sup_{0 \leq t \leq T} \int_{\vert x \vert \geq R} \vert x \vert^p \, d[S_hu^0]_x(t,x) 
\leq \varepsilon
$$
by Proposition \ref{equiintegrable}. Let then $ \varphi \in {\cal C}^{\infty}_{c}(\R)$ such that 
$0 \leq \varphi \leq 1$ and $\varphi(x) =0$ if $\vert x \vert \leq R$.

On one hand
$$
\int_{\R} \varphi(x) \vert x \vert^p \, d[S_hu^0]_x(t,x) \,
 \to \int_{\R} \varphi(x) \vert x \vert^p \, du_x(t,x)
 $$
 as $h$ goes to $0$ since $\varphi(x) \vert x \vert^p \in {\cal C}^{\infty}_{c}(\R)$ and $[S_hu^0]_x(t,.)$ tends 
 to $u_x(t,.)$ in distribution sense.
 On the other hand
$$
\int_{\R} \varphi(x) \vert x \vert^p \, d[S_hu^0]_x(t,x) \, \leq \, \varepsilon
$$
for all $0 \leq h,t \leq T$. Hence at the limit
$$
\int_{\R} \varphi(x) \vert x \vert^p \, du_x(t,x) \, \leq \, \varepsilon
$$
for all $t \leq T$, from which it follows that
$$
\sup_{0 \leq t \leq T} \int_{\vert x \vert \geq R} \vert x \vert^p \, du_x(t,x) \, \leq \, \varepsilon,
$$
which means that indeed $\ds \sup_{0 \leq t \leq T} \int_{\vert x \vert \geq R} \vert x \vert^p \, du_x(t,x)$ goes to $0$
as $R$ goes to infinity.
 
\medskip
 
 From these two results we deduce the continuity result by Proposition \ref{cvWp/weak}.

 \subsection{Proof of Corollary \ref{corsolutionU1}}
 \indent
 
 Given $t \geq 0$, $S_hu^0(t,.)$ converges to $u(t,.)$ in $L^{1}_{loc}(\R)$ by Proposition \ref{schemeconvergence},
 so for all $s,t, n \geq 0$ we have
$$
\Vert u(t,.) - u(s,.) \Vert_{L^1([-n,n])} = \lim_{h \to 0} \Vert S_hu^0(t,.) - S_hu^0(s,.) \Vert_{L^1([-n,n])}.
$$

But
$$
\Vert S_hu^0(t,.) - S_hu^0(s,.) \Vert_{L^1([-n,n])} \leq \Vert S_hu^0(t,.) - S_hu^0(s,.) \Vert_{L^1(\R)}
\leq \vert t-s \vert \, \Vert f' \Vert_{L^{\infty}(\R)}
$$
for all $h \geq 0$ by Proposition \ref{contractL1}, so letting $h$ go to $0$ we get
$$
\Vert u(t,.) - u(s,.) \Vert_{L^1([-n,n])} \leq \vert t-s \vert \, \Vert f' \Vert_{L^{\infty}(\R)}.
$$
Since this holds for all $n \geq 0$, we obtain Corollary \ref{corsolutionU1}.

\section{Extension to viscous conservation laws}\label{sectionviscous}
\indent

In this section we let $\nu$ be a positive number and consider the viscous conservation law
\begin{equation}\label{viscouscl}
u_t + f(u)_x =  \nu \, u_{xx} \qquad t> 0, \; x \in \R 
\end{equation}
with initial datum $u^0 \in L^{\infty}(\R)$.

Assuming that $f$ is a locally Lipschitz real-valued function on $\R$, and calling 
{\it solution} 
a function $u$ in $L^{\infty}([0,+\infty[ \times \R)$ such that \eqref{viscouscl} holds
in the sense of
distributions, it is known that, given $u^0 \in L^{\infty}(\R)$, there exists a unique 
solution $u$ to \eqref{viscouscl}. If 
moreover $u^0 \in \cal U$, then $u(t,.)$ also belongs to $\cal U$
for all $t \geq 0$, and the $W_p$ contraction property stated in Theorem \ref{theoreme}
in the inviscid case $\nu = 0$ still holds:
\begin{theo}\label{theoremevisqueux}
Given  a locally Lipschitz real-valued function $f$ on $\R$ and two initial data $u^0$ and ${\tilde u}^0$ 
in ${\cal U}$, let $u$ and 
${\tilde u}$ be the associated solutions to \eqref{viscouscl}. Then, for any 
$t \geq 0$ and $p \geq 1$, we have (with possibly infinite values)
$$
W_p(u_{x}(t,.), {\tilde u}_{x}(t,.)) \, \leq \, W_p(u_{x}^{0}, {\tilde u}_{x}^{0}).
$$
\end{theo}

\medskip

We briefly mention how this contraction property for the viscous conservation law allows to
recover the same property for the inviscid equation, given in Theorem \ref{theoreme}. Given 
some initial datum $u^0$ in $\cal U$ and $\nu >0$, let indeed $u_{\nu}$ be the corresponding
solution to the viscous equation \eqref{viscouscl}. Then it is known (see \cite{S} for instance)
that $u_{\nu}(t,.)$ converges in $L^{1}_{loc}(\R)$ to the solution $u(t,.)$ to the
inviscid conservation law \eqref{cons law} with initial datum $u^0$. From this the argument
already used in Section \ref{soussectionpreuve} (with $S_hu^0(t,.)$ intead of $u_{\nu}(t,.)$)
enables to recover Theorem \ref{theoreme}.

\medskip

The proof of Theorem \ref{theoremevisqueux} follows the lines of Sections 
\ref{sectiondiscretisation} and \ref{sectionpreuve}
and makes use of a time-discretization of equation \eqref{viscouscl} based on the discretization
of the inviscid conservation law
previously discussed. More precisely, given a time step $h > 0$, we first
map $u^0 \in {\cal U}$ to $T_hu^0$ as in section 4.1, and then let $T_hu^0$ evolve along the heat
equation on a time
interval $h$, that is, map it to 
$$
\T_hu^0 = K_h * T_hu^0 
$$
where $K_h$ is the heat kernel defined  on $\R$ by
$$
K_h(z) = \frac{1}{\sqrt{4 \pi h}} \, e^{-\frac{z^2}{4h}}.
$$
Then, defining an approximate solution $\cal S_hu^0$ by iterating the $\T_h$ operator as in Definition
\ref{defiSh}, we prove that Propositions \ref{propTh} and \ref{Shcontract} still hold for the
new $\T_h$ and $\cal S_h$ operators (Note that the convolution with the heat kernel is a contraction for the Wasserstein distance of any finite order).

Proposition \ref{contractL1} only holds
assuming that $v$ is twice derivable with $v''$ in $L^1(\R)$: it more precisely reads
$$
\Vert \cal S_hv(t,.) - \cal S_hv(s,.) \Vert_{L^1(\R)} \leq \, 
\vert t-s \vert \, [\Vert f' \Vert_{L^{\infty}(]0,1[)}  +  \Vert v'' \Vert_{L^1(\R)}].
$$
As in Section \ref{cvL1loc}, this enables to prove that, given $u^0$ in $\U$, twice derivable with
$(u^0)''$ in $L^1(\R)$, the family $(\cal S_hu^0)_h$ converges in $\cal C([0,+\infty[,L^{1}_{loc}(\R))$
to the solution of \eqref{viscouscl} with initial datum $u^0$.

With this convergence result in hand we follow the lines of Section \ref{soussectionpreuve} to prove
Theorem \ref{theoremevisqueux} in the case of twice derivable initial data, with $L^1$ second derivative,
while the general case follows by a density argument.

\bigskip

\bigskip

\noindent
{\bf Acknowledgments:}
This work has been partly supported by the European network Hyke, funded by the EC contract 
HPRN-CT-2002-00282.


\medskip

\begin{thebibliography}{99}

\bibitem[1]{B} Y. Brenier, R\'esolution d'\'equations d'\'evolution quasilin\'eaires en
                dimension {$N$} d'espace \`a l'aide d'\'equations lin\'eaires en dimension {$N+1$},
                 {\it J. Diff. Eq.} {\bf 50}, 3, p. 375-390 (1983).


\bibitem[2]{BG} Y. Brenier and E.  Grenier, Sticky particles and scalar conservation laws, {\it SIAM J. Numer. Anal.}, {\bf 35}, 6, p. 2317-2328 (1998).

     
\bibitem[3]{H} L. H{\"o}rmander, {\it Lectures on  nonlinear hyperbolic differential equations}, 
Mathematics \&  Applications (26), Springer, Berlin (1997).

\bibitem[4]{K} S. Kru$\mathrm {\check z}$kov, Generalized solutions of the Cauchy problem in the large for
  first order nonlinear equations, {\it Dokl. Akad. Nauk. SSSR}, {\bf 187}, p. 29-32 (1969).

\bibitem[5]{S} D. Serre, {\it Syst\`emes de lois de conservations I}, Diderot, Paris (1996).

\bibitem[6]{V} C. Villani, {\it Topics in optimal transportation}, Grad. Stud. Math (58), American Mathematical
Society, Providence (2003).

\end{thebibliography}
\end{document}